\documentclass{article}

\usepackage[preprint]{neurips_2025}

\usepackage{times}  % DO NOT CHANGE THIS
\usepackage{helvet}  % DO NOT CHANGE THIS
\usepackage{courier}  % DO NOT CHANGE THIS
\usepackage{graphicx} % DO NOT CHANGE THIS
\usepackage{natbib}  % DO NOT CHANGE THIS AND DO NOT ADD ANY OPTIONS TO IT
\usepackage{caption} % DO NOT CHANGE THIS AND DO NOT ADD ANY OPTIONS TO IT
\usepackage{algorithm}
\usepackage{algorithmic}

\usepackage{todonotes} % for comments

\usepackage{booktabs} 
\usepackage{amsfonts}   
\usepackage{amsmath} 
\usepackage{nicefrac}     
\usepackage{microtype}     
\usepackage{mdframed}

\usepackage{scalerel}
\usepackage{txfonts}
\usepackage{fontawesome}
\usepackage{todonotes}
\usepackage{xcolor}
\usepackage{multirow}
\usepackage{bm}
\usepackage{paralist}
\usepackage{enumitem}
\usepackage{caption}[small]
\usepackage{subcaption}

\usepackage[utf8]{inputenc} 
\usepackage[T1]{fontenc}   
\usepackage{url}           
\usepackage{booktabs}       
\usepackage{amsfonts}    
\usepackage{nicefrac}     
\usepackage{microtype}    
\usepackage{xcolor}      
\usepackage{wrapfig}

\DeclareMathOperator*{\argmin}{argmin}

 \newcommand{\cB}{\mathcal{B}} 
\newcommand{\cC}{\mathcal{C}} 
 \newcommand{\cF}{\mathcal{F}}
\newcommand{\cG}{\mathcal{G}} 
 \newcommand{\cL}{\mathcal{L}}
 \newcommand{\cN}{\mathcal{N}}
\newcommand{\cO}{\mathcal{O}} \newcommand{\cP}{\mathcal{P}}
 \newcommand{\cS}{\mathcal{S}}

\usepackage{newfloat}
\usepackage{listings}
\DeclareCaptionStyle{ruled}{labelfont=normalfont,labelsep=colon,strut=off} % DO NOT CHANGE THIS
\floatstyle{ruled}
\newfloat{listing}{tb}{lst}{}
\floatname{listing}{Listing}
\title{Learning to Solve Constrained Bilevel Control Co-Design Problems}

\author{%
  James Kotary \\
  Pacific Northwest National Laboratory\\
  \texttt{james.kotary@pnnl.gov} \\
  \And
  Himanshu Sharma \\
  Pacific Northwest National Laboratory\\
\texttt{himanshu.sharma@pnnl.gov} \\
  \AND
  Ethan King \\
  Pacific Northwest National Laboratory \\
  \texttt{ethan.king@pnnl.gov} \\
  \And
  Draguna Vrabie\\
  Pacific Northwest National Laboratory \\
\texttt{draguna.vrabie@pnnl.gov} \\
  \And
  Ferdinando Fioretto\\
  University of Virginia\\
  \texttt{fioretto@virginia.edu} \\
  \And
  Jan Drgona\\
  Johns Hopkins University\\
  \texttt{jdrgona1@jh.edu} \\
}

\begin{document}

\maketitle\allowdisplaybreaks\sloppy

\begin{abstract}
We propose a learning to optimize (L2O) method for solving  \emph{constrained parametric bilevel problems} that arise in control co-design, where upper-level design variables are coupled with lower-level optimal control through explicit \emph{coupling constraints}. Our self-supervised framework comprises: (i) a \emph{differentiable optimization layer} to enforce lower-level optimality, and (ii) a \emph{differentiable gradient-based projection} routine that iteratively reduces coupling-constraint violation while maintaining feasibility of upper-level constraints. A soft penalty is used during training to initialize predictions near feasibility, enabling stable end-to-end learning.
On bilevel QPs with certified optima, our learned models achieve $10^{-2}$ relative optimality gaps while running $\sim 10^{2}\times$ faster than a mixed-integer programming (MIP) reformulation. On two optimal control co-design tasks, our approach yields $15$–$19\%$ lower design cost and $\sim 10^{4}\times$ faster inference
than a particle swarm optimization (PSO) baseline, while maintaining comparable constraint satisfaction. 
These results indicate that the proposed L2O method can deliver real-time, high-quality approximations for challenging bilevel programming problems that are computationally prohibitive using conventional methods.
\end{abstract}

% \begin{abstract}
% Recent work in machine learning has shown great promise in deep neural networks as fast approximators of difficult optimization problems. Such advances are key to providing real-time, high-quality solutions increasingly demanded in many engineering applications. Prior Learning-to-Optimize (L2O) methods focus almost entirely on single-level programs, in contrast to bilevel programs, whose constraints are themselves expressed in terms of optimization subproblems. Bilevel programs have numerous important use cases but are notoriously difficult to solve, particularly under stringent time demands. This paper proposes a framework for learning to solve a class of challenging bilevel optimization problems by leveraging modern techniques for differentiation through optimization problems. The framework is illustrated on an array of synthetic bilevel programs, as well as challenging control system co-design problems, showing how neural networks can be trained as efficient approximators of parametric bilevel optimization.

% \end{abstract}

\section{Introduction}
\label{sec:Introduction}

Bilevel optimization problems arise in a wide range of applications, from economics and game theory \citep{binmore2007playing}, to operations management and logistics \citep{sadigh2012manufacturer}, as well as engineering system design \citep{bergonti2024co}. Unfortunately, bilevel problems are, in general, notoriously difficult to solve. They are typically NP-hard, and depending on their particular qualities, may even lack efficient frameworks for their approximate solution \citep{beck2021gentle}. This complexity challenges their use in applications requiring \emph{real-time or repeated solutions}. Yet, many applications of bilevel optimization demand repeated solution of related problem instances, such as when engineering design decisions are considered across a variety of scenarios and objectives. 

This paper introduces a novel framework for applying deep learning to solve a broad class of bilevel problems, which include  \emph{coupling constraints}, i.e., constraints that bind upper and lower level decision variables. The approach uses differentiable optimization solvers at the problem's lower level, allowing for gradient-based training of neural networks to approximate the upper-level solution, as well as internal  \emph{correction} routines to enforce coupling constraints between the upper and lower-level problems. The resulting models are trained to map the parameters of a bilevel problem to its optimal solution. We motivate this approach primarily using problems from optimal control co-design: a setting which calls for the optimization of engineering systems, with respect to system design objectives at the upper level and subject to system dynamics determined by optimal control problems at the lower level. We demonstrate how the proposed framework can be used to learn solutions in response to varying design objectives and desiderata. 

\noindent \textbf{Contributions.} The main contributions of this paper are:
\textbf{(1)} We propose a novel Learning to Optimize method for bilevel optimization problems, whose core concept is based on using differentiable optimization to ensure optimality of the lower-level solution, while ensuring constraint satisfaction at the upper level via differentiable projections.
\textbf{(2)}  We show how differentiable optimization can also be used to compose gradient-based correction mechanisms for the satisfaction of \emph{coupling constraints} - a significant challenge in bilevel problems.
\textbf{(3)} We demonstrate the proposed framework on a collection of challenging constrained bilevel optimization and control-co-design problems, including those with \emph{nonconvex} lower level problems and complex constraints. In particular, we demonstrate its ability to learn high-quality approximate parametric solutions for a distribution of problems. Full code is found at: \url{https://github.com/anon-research-coding/Control-Design-Learning}
\vspace{-8pt}
\section{Related Work}
\label{sec:RelatedWork}
%This work is at the intersection of bilevel optimization and learning to optimize. Related work within each area is reviewed below. \rev{do we need this sentence pre-subsections or can we save the space?}
\noindent \textbf{Bilevel optimization methods.}
Bilevel optimization problems are notoriously difficult to solve, and in general have no efficient solution methods except in special cases. Particularly when problems are sufficiently small and a convex structure is present, generic solution methods tend to focus on single-level reformulations (via the KKT conditions or the optimal value function), which can be solved by mixed-integer programming \citep{beck2021gentle}. When a special structure is present, such as linearity at both levels or a lack of coupling constraints, gradient methods have been proposed that tend to rely on penalty methods \citep{cerulli2021solving, ghadimi2018approximation, solodov2007explicit}. 
Recently, \citep{sharifi2025safegradientflowbilevel,abolfazli2025} solved a special class of bilevel problems using upper-level gradient descent with a quadratic programming (QP)-based safety filter for enforcing the KKT conditions of the lower-level problem.
The complex structure possessed by more general bilevel problems has led to greater interest in metaheuristics, such as evolutionary algorithms \citep{sinha2013efficient,bergonti2023co} and particle swarm optimization \citep{li2006hierarchical} over classical alternatives. A comprehensive review of modern approaches is found in \cite{sinha2017review}.

\noindent \textbf{Learning to Optimize.} L2O is a subfield of ML concerned with learning \emph{fast approximations} to challenging optimization problems~\citep{vanhentenryck2025}. One distinct branch focuses on learning information to accelerate the convergence of a classical solver \citep{bengio2020machine}. In continuous optimization, these include prediction of active constraints \citep{ng2018statistical}, optimization problem parameters~\citep{Agrawal2021}, stepsizes \citep{amos2023tutorial}, primal variables \citep{JMLR:v25:23-1174,Bertsimas2022}, and noneuclidean metrics \citep{king2024metric}, while for integer variables they include branching rules \citep{khalilbranch}, cutting planes \citep{tang2020reinforcement}, variable partitions in large neighborhood search \citep{SongLNS}, or primal variables~\citep{tang2025,Huang2025}. A separate branch aims at training ML models to predict optimal solutions directly from a representation of the problem \citep{kotary2021end}. 

A key challenge for such an approach is to maintain the feasibility of the predicted solutions to arbitrary constraints. Proposals for addressing this aspect are based on differentiable projections \citep{wilder2018melding}, reparametrization tricks \citep{frerix2020homogeneous,konstantinov2024new}, dual-variable estimation \citep{fioretto2020lagrangian,park2023self}, and gradient-based constraint correction routines \citep{donti2021dc3}. Finally, some recent works have proposed to accelerate bilevel optimization with ML. \citep{dumouchelle2024neur2bilo} proposes in the non-parametric case to learn an optimal value reformulation from solved examples with a ReLU network, which is embedded into an MIP solver. \citep{shen2020learning,Andrychowicz2016} propose to learn solutions of a parametric bilevel program directly, albeit without any constraints at either level. \emph{This paper extends the toolkit for direct learning of solutions to the case of bilevel optimization with continuous variables and a full set of constraints, including coupling constraints}.
%\jan{Lets put together conceptual diagram of the whole methodology. I can share my svg files from previous papers to adapt. }

\section{Problem Setting}
\label{sec:Setting}
The goal of this paper is to learn to solve parametric bilevel optimization problems represented by \eqref{UL}. We take the convention that optimization variables $\bm{y},\bm{z}$ are written as function arguments while problem parameters $\bm{p}$ are written as subscripts. We consider a pair of problems: 
\vspace{-4pt}
\begin{figure}[h!]
\centering
\begin{minipage}{0.48\textwidth}
\vspace{-\baselineskip}
\begin{subequations}\label{UL}\begin{align}
         \cB(\bm{p}) \coloneqq \underset{\bm{y}}{\argmin} \;\;\;\; & {\cL}_{\bm{p}}(\bm{y},\bm{z}) \label{UL-obj} \\
         \textit{s.t.}\;\;\;\;
         & \bm{z} \in \cO_{\bm{p}}(\bm{y}) \label{UL-opt} \\
         & \bm{y} \in \mathcal{C}_{\bm{p}} \label{UL-constr} \\
         & \bm{U}_{\bm{p}}(\bm{y},\bm{z}) \leq \bm{0} \label{UL-coupling}
    \end{align}\end{subequations}
\vspace{-\baselineskip}
\end{minipage}
\hfill
\begin{minipage}{0.48\textwidth}
\vspace{-\baselineskip}
\begin{subequations}
    \label{LL}
    \begin{align}
         \cO_{\bm{p}}(\bm{y}) \coloneqq \underset{\bm{z}}{\argmin} \;\;\;\; & {l}_{\bm{p}}(\bm{y},\bm{z}) \label{LL-obj} \\
         \textit{s.t.}\;\;\;\;
         & \bm{z} \in \mathcal{S}_{\bm{p}}(\bm{y}). \label{LL-constr}
    \end{align}\end{subequations}
\end{minipage}
\vspace{-\baselineskip}
\end{figure}

%\begin{subequations}
%    \label{UL}
%    \begin{align}
%         \cB(\bm{p}) = {\argmin}_{\bm{y}} \;\;\;\; & {\cL}_{\bm{p}}
%         (\;\bm{y},\bm{x},\bm{u}\;)\label{UL-obj} \\
%         \textit{s.t.}\;\;\;\; & \bm{y} \in \mathcal{C} \label{UL-constr} \\
%         % & \bm{h}(\;\bm{p},\bm{y}\;) = \bm{0} \label{UL:params} \\
%         & (\;\bm{x},\bm{u}\;) \in \cO(\bm{p},\bm{y}) \label{UL-opt} \\
%                  & U(\bm{x},\bm{y}) \leq \bm{0} \label{ULC-coupling}.
%    \end{align}
%\end{subequations}
\noindent The defining feature of problem \eqref{UL} is its constraint \eqref{UL-opt}, which holds that the variables $\bm{z}$ be the solution to another optimization problem \eqref{LL}, which in turn depends on the variables $\bm{y}$. We call problem \eqref{UL} the \emph{upper-level problem}, and \eqref{LL} is the \emph{lower-level problem}, while $\bm{y}$ and $\bm{z}$ are the \emph{upper} and \emph{lower-level} variables, respectively. Our goal is to learn a fast approximator that solves the coupled problems (\ref{UL},\ref{LL}), over a distribution of problem parameters denoted as $\mathcal{P}$.
%Our goal is to \emph{learn the mapping} $\bm{p} \to \cB(\bm{p})$ with a neural network, over a distribution  of problem parameters $\bm{p} \sim \mathcal{P}$.

\noindent \textbf{Classes of constraints.} We distinguish three sets of constraints at the upper level. The condition \eqref{UL-constr} constrains only the upper-level variables, while \eqref{UL-opt} prescribes $\bm{z}$ as a solution to the lower-level problem \eqref{LL} resulting from $\bm{y}$. Additionally, the \emph{coupling constraints} \eqref{UL-coupling} significantly complicate the solution of \eqref{UL}. They impose additional conditions on the relationship between upper and lower-level variables in \eqref{UL}, preventing solution concepts based on their separation or decoupling \citep{beck2021gentle}. A large portion of algorithms for bilevel programming cannot accommodate problems with coupling constraints \citep{sinha2017review}, and to the best of our knowledge, \emph{no previous work has ventured to propose an L2O framework for learning to solve them parametrically}.
%The \emph{upper-level objective}  ${\cL}_{\bm{p}}$ is a joint function of both $\bm{y}$ and its associated $\bm{z}$.

 %The next section details the Optimal Control Co-design setting as our motivating application. On the other hand, we emphasize that the paper's proposed methodology is \emph{generic} and applies to a broad class of bilevel programs. We do however make one crucial assumption throughout the paper: that 

%\rev{make bulleted points about main challenges earlier, including infeasibility at lower level (Section 4.4) - even in intro-contributions}

\noindent \textbf{Conditions on the problem form.} The proposed framework aims to learn solutions to a broad class of problems having the form (\ref{UL},\ref{LL}). However, it depends on a key condition: for all $\bm{y} \in \mathcal{C}_{\bm{p}}$ and $\bm{p} \sim \mathcal{P}$, the solution to the lower-level problem \eqref{LL} must be  \emph{unique} whenever it \emph{exists} - existence is not required but assumed to hold without loss of generality until it is relaxed using a reformulation trick, introduced in Section \ref{sec:LowerLevelConstraints}. This leads $\cO_{\bm{p}}(\bm{y})$ to define a function, which we further suppose to be \emph{differentiable}. These conditions are not met in all bilevel problems. However, they are satisfied in many important cases, such as when \eqref{LL} is a \emph{Model Predictive Control} (MPC) problem. Several recent works have proposed efficient techniques for differentiating through MPC  \citep{amos2018differentiable,dinev2022differentiable}, which generally have unique solutions given sufficiently specified objectives.
The functions $\cL_{\bm{p}}$ and $\bm{U}_{\bm{p}}$ should also be differentiable. Importantly, though, a.e. or pseudo-differentiability is often sufficient for effective gradient descent optimization, as in the case of ReLU activations. Following common practice in L2O and machine learning, these conditions are treated as practical guidelines for applying the proposed architecture, rather than formal requirements.
%\jan{there is a lot of white space in eq. 3 - we can wrap it around the text}
\section{Learning to Solve Bilevel  Optimization}
This section presents a self-supervised method for learning  to solve the parametric bilevel optimization problems described in Section \ref{sec:Setting}. In particular, it trains an ML model to approximate the mapping \eqref{UL}, from problem parameters $\bm{p}$ to \emph{upper-level} solutions $\bm{y}^{\star} = \cB(\bm{p})$. Let $\hat{\cB}_{\theta}$ denote that hypothetical model, with weights $\theta$. Assume a training set  of $n_p$ problem instances, each specified by a vector of problem parameters $\{\bm{p}_{(i)}\}_{i=1}^{n_p}$ which are drawn from the distribution $\cP$. A training procedure for $\hat{\cB}_{\theta}$ should minimize the objective function $\cL_{\bm{p}}(\hat{\bm{y}})$ attained by its predicted solutions  $\hat{\bm{y}} = \hat{\cB}_{\theta}(\bm{p})$ in expectation over $\cP$, resulting in the empirical risk minimization (ERM)~\eqref{eq:ERM}.
\iffalse
\begin{subequations}
\label{eq:ERM}
\begin{align}
    \min_{\bm{\theta}}  \;\; \underset{\bm{p} \sim \cP}{\mathbb{E}}\;\; \Big[ & \cL_{\bm{p}}(\hat{\bm{y}}, \hat{\bm{z}}) \Big] \label{eq:ERM_objective} \\
    \textit{s.t.} \;\;\;
    & \hat{\bm{z}} \in \cO_{\bm{p}}(\bm{y})   \;\;\;\;\;\;\;\;\;\;\; \forall \bm{p} \sim \cP 
     \label{eq:ERM_c2} \\
    & \hat{\bm{y}} \in \mathcal{C}_{\bm{p}} \;\;\;\;\;\;\;\;\;\;\;\;\;\;\; \forall \bm{p} \sim \cP \label{eq:ERM_c1}\\
    & \bm{U}_{\bm{p}}(\hat{\bm{y}},\hat{\bm{z}}) \leq \bm{0} \;\;\;\;\;\;\; \forall \bm{p} \sim \cP, \label{eq:ERM_c3}
\end{align}
\end{subequations}
\fi
\begin{wrapfigure}{r}{0.40\textwidth} % 'r' = right; adjust width if needed
\vspace{-2em} % tighten space above the equation
\begin{subequations}
\label{eq:ERM}
\begin{align}
\min_{\bm{\theta}} \; 
    \mathbb{E}_{\bm{p} \sim \cP}\!\left[
        \cL_{\bm{p}}(\hat{\bm{y}}, \hat{\bm{z}})
    \right] 
    & \label{eq:ERM_objective} \\[-0.3em]
\text{s.t.}\quad 
    \hat{\bm{z}} \in \cO_{\bm{p}}(\hat{\bm{y}}) 
    & \quad \forall \bm{p} \sim \cP, 
    \label{eq:ERM_c2} \\
    \hat{\bm{y}} \in \mathcal{C}_{\bm{p}} 
    & \quad \forall \bm{p} \sim \cP, 
    \label{eq:ERM_c1} \\
    \bm{U}_{\bm{p}}(\hat{\bm{y}}, \hat{\bm{z}}) \le \bm{0} 
    & \quad \forall \bm{p} \sim \cP.
    \label{eq:ERM_c3}
\end{align}
\end{subequations}
\vspace{-2em} % tighten space below
\end{wrapfigure}
Constraints (\ref{eq:ERM_c1},\;\ref{eq:ERM_c2},\;\ref{eq:ERM_c3}) require that each such output and its resulting pair $(\hat{\bm{y}}, \hat{\bm{z}})$ is a feasible solution to the bilevel problem (\ref{UL},\ref{LL}), where $\hat{\bm{y}} \coloneqq \hat{\cB}_{\theta}(\bm{p})$. 
Learning solutions subject to such complex constraints is inherently challenging, as any predicted $\hat{\bm{y}}$ and its corresponding $\hat{\bm{z}}$ are unlikely to satisfy the coupling constraint \eqref{eq:ERM_c3} after solving \eqref{eq:ERM_c2}, even if $\hat{\bm{y}} \in \mathcal{C}_{\bm{p}}$ as required in \eqref{eq:ERM_c1}. In this section, we propose an architecture for the model $\hat{\cB}_{\theta}$, along with a method for training it to approximate the ERM goal \eqref{eq:ERM}.
%Its core concept is based on identifying optimization subproblems which, when treated as \emph{differentiable functions}, may be used to compose an end-to-end trainable predictor of feasible solution pairs $(\hat{\bm{y}},\hat{\bm{z}})$. Specifically, it adopts a differentiable subroutine \redtext{use this to mention DC3/ correction or just remove} which iteratively refines predicted solutions to satisfy the coupling constraint \ref{eq:ERM_c3} while \emph{maintaining feasibility} to the constraints \eqref{eq:ERM_c2} and \eqref{eq:ERM_c1} at each step. Before prescribing the full training method, the main architectural components of $\hat{\cB}_{\theta}$ are introduced next.
\subsection{Satisfying Constraints with Differentiable Optimization Modules}
The proposed architecture employs differentiable optimization to implement two of its main components, prefaced below. In general, the derivatives between a problem's parameters and its optimal solution can be found by implicit differentiation, either of the KKT conditions \citep{amos2019optnet, gould2021deep,agrawal2019differentiable, agrawal2019differentiating} or a fixed-point condition \citep{JMLR:v25:23-1174}. This proposal is agnostic as to which implementation is used. %The purpose of differentiable optimization is often to use optimization problems as modules in trainable machine learning models, and we employ two such modules within $\hat{\cB_{\theta}}$.

\noindent \textbf{Differentiable solution of the lower-level problem.}
For a given $\bm{p}$ and any predicted upper-level solution $\hat{\bm{y}}$, a differentiable solver of problem \eqref{LL}  produces $\hat{\bm{z}} = \cO_{\bm{p}}(\hat{\bm{y}})$ along with $ \frac{\partial \hat{\bm{z}}}{\partial \hat{\bm{y}}}$. This provides a means by which feasible solution pairs \eqref{eq:ERM_c2} can be computed and back-propagated as part of an end-to-end trainable model.
In the applications of Section \ref{sec:Experiments}, it corresponds to generating an optimal control policy $\hat{\bm{z}}$ as a differentiable function of learned design parameters $\hat{\bm{z}}$.

%In particular, this section presents a differentiable model for learning feasible pairs $(\hat{\bm{y}},\hat{\bm{z}})$, which iteratively refines the initial prediction of a neural network. Importantly, it maintains feasibility to the constraints \eqref{eq:ERM_c2} and \eqref{eq:ERM_c1} at each step, while iterating toward feasibility to the more difficult coupling constraint \eqref{eq:ERM_c3}.
%\rev{commented paragraph here that should be checked}
%\rev{left off here}
%Given any predicted upper-level solution, we may restore its feasibility to \eqref{UL-constr} by a differentiable projection onto $\mathcal{C}_{\bm{p}}$. Subsequently, lower-level solutions can be made to satisfy the the condition \eqref{UL-opt} by differentiably solving problem \eqref{LL-discr}. Differentiability of each of these elements allows the unsupervised loss function ${\cL}_{\bm{p}}$ \rev{write the ERM above} to be minimized by end-to-end gradient descent training while maintaining feasibility. Section \ref{sec:LearningUpperLevel} introduces this idea by applying it to simplified case of problem \eqref{UL} which is absent of coupling constraints \eqref{UL-coupling}.  Section \ref{sec:Coupling_Constraints} then builds on this concept to incorporate coupling constraints \eqref{UL-coupling}, by composing those same elements into a differentiable \emph{coupling constraint correction} routine.
%\subsection{Warm-up: Learning Bilevel Optimization without Coupling}
\noindent \textbf{Differentiable projection at the upper level.}
In addition to satisfying constraint \eqref{eq:ERM_c2} with its solution pair, any candidate $\hat{\bm{y}}$ must also satisfy its own upper-level constraint \eqref{eq:ERM_c1}. This constraint is handled using a projection operator $\bm{\Pi}_{\mathcal{C}}$.
The defining problem \eqref{C_projection} can be solved in a differentiable library such as \textit{cvxpylayers} \citep{agrawal2019differentiable}, but automatic differentiation in PyTorch suffices when a closed form is available (e.g., the projection onto $\{ \bm{x}: \bm{x} \geq 0 \}$ is well-known to be a ReLU function). 
\begin{wrapfigure}{r}{0.3\textwidth} 
\vspace{-1.5em} 
\begin{equation}
\label{C_projection}
    \bm{\Pi}_{\mathcal{C}}(\bm{y}) 
    = \argmin_{\bm{w} \in \mathcal{C}} 
      \| \bm{y} - \bm{w} \|_2^2
\end{equation}
\vspace{-2em} % reduce space below the equation
\end{wrapfigure}
Differentiable projections are a mainstay tool for constraint satisfaction in L2O \citep{sambharya2023end, king2024metric, wilder2018melding}. Besides projections, other mechanisms could be used to guarantee constraint satisfaction \citep{chen2023end,tordesillas2023rayen}. However, the projection operator is chosen for its role in the algorithm presented next, resulting in the \emph{projected gradient descent} method, which has well-studied convergence properties \citep{beck2017first}.

\subsection{End-to-End Trainable Architecture}
\label{sec:Coupling_Constraints}
\begin{wrapfigure}{R}{0.5\textwidth}
    \vspace{-6pt}
    \centering
    \includegraphics[width=0.95\linewidth]{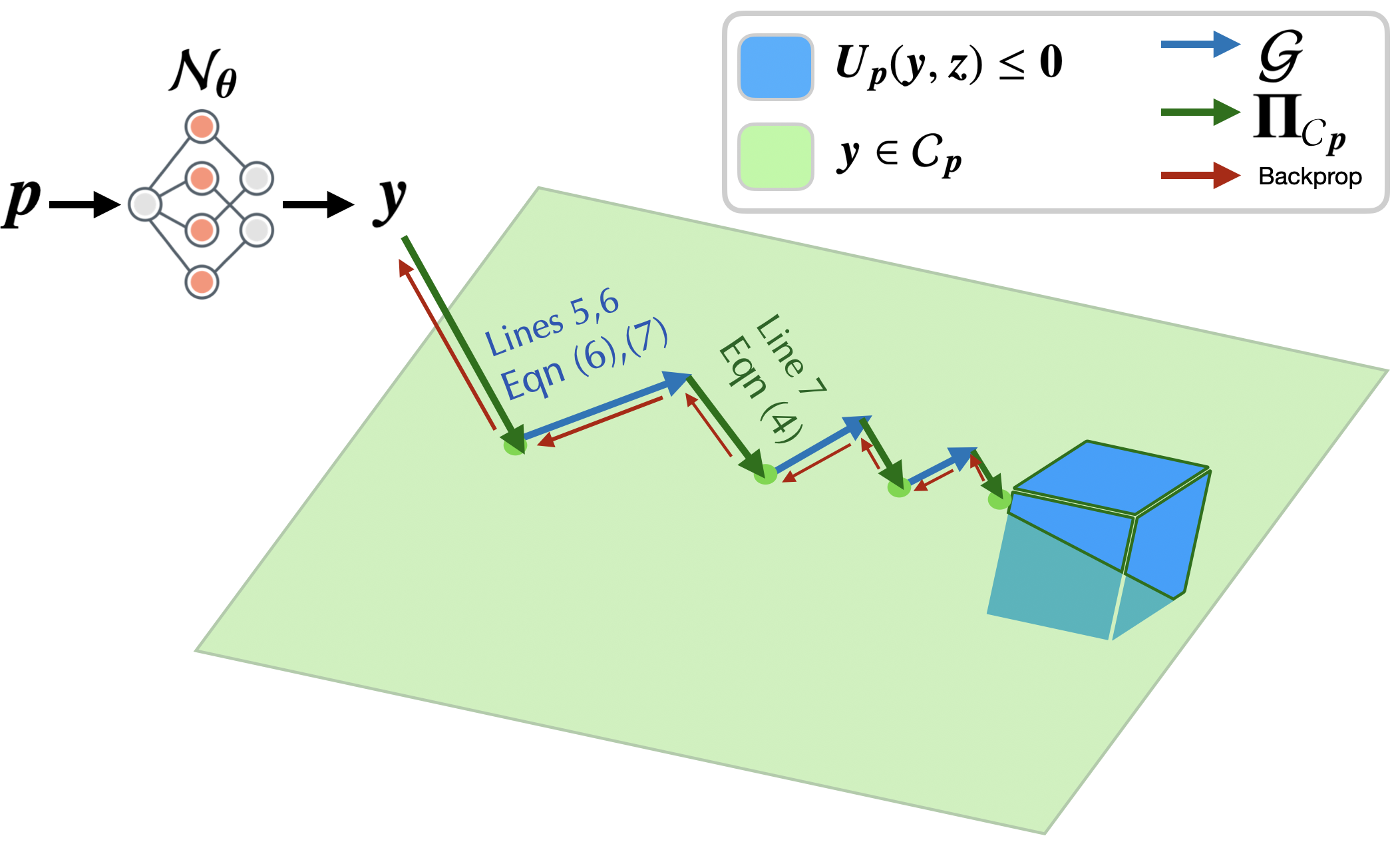}
    \caption{Illustrating the end-to-end trainable model with algorithm line and equation numbers. }
    \label{fig:illustration}
    \vspace{-1pt}
\end{wrapfigure}

We can now present the complete architecture and training routine of a model $\hat{\cB}_{\theta}$ which learns to solve problem \eqref{eq:ERM}. Let $\cN_{\theta}$ be a deep neural network with weights $\theta$, which predicts initial estimates $\hat{\bm{y}} = \cN_{\theta}(\bm{p})$  of an upper-level solution. Composition of $\cN_{\theta}$ with $\bm{\Pi}_{\mathcal{C}_{\bm{p}}}$ ensures feasibility to \eqref{eq:ERM_c1}, and further composition with $\cO_{\bm{p}}$ produces a solution pair $( \hat{\bm{y}},\cO_{\bm{p}}(\hat{\bm{y}})  ) $ which satisfies \eqref{eq:ERM_c2}. Therefore the function  $\cO_{\bm{p}} \circ \bm{\Pi}_{\mathcal{C}_{\bm{p}}}$
can be viewed as one which maps infeasible upper-level estimates to solution pairs $(\hat{\bm{y}},\hat{\bm{z}})$ satisfying  \eqref{eq:ERM_c2} and  \eqref{eq:ERM_c1}, but not the coupling constraint $\bm{U}_{\bm{p}}(\hat{\bm{y}},\hat{\bm{z}}) \leq \bm{0}$ in \eqref{eq:ERM_c3}. We define the \emph{Coupling Constraint Violation} as in \eqref{eq:CCV}, along with the gradient of its Euclidean norm in \eqref{eq:violation_grad}. Note that  ${\bm{y}}$ is recognized as the independent variable and ${\bm{z}} = \cO_{\bm{p}}({\bm{y}})$ as dependent. 
\begin{minipage}{0.36\textwidth}
\vspace{-0.6\baselineskip}
\begin{equation}
    \label{eq:CCV}
       \nu({\bm{y}}) \coloneqq  \textsc{ReLU}\left( \bm{U}(\; {\bm{y}},\cO_{\bm{p}}({\bm{y}})  \;) \right),
\end{equation}
\end{minipage}
\hfill
\begin{minipage}{0.62\textwidth}
\vspace{-0.6\baselineskip}
\begin{align}
    \label{eq:violation_grad}
       \nabla_{\bm{y}} \| \nu(\bm{y}) \|^2   
        = 2 \nu(\bm{y})\frac{d \nu}{d\bm{y}} 
        = 2 \nu(\bm{y}) \frac{d\nu}{d\bm{U}}  \left[   \frac{\partial \bm{U}}{\partial c} + \frac{\partial \bm{U}}{\partial \mathcal{O}_{\bm{p}}}  \frac{\partial  \mathcal{O}_{\bm{p}}}{\partial \bm{y}}  \right].
\end{align}
\end{minipage}
The nontrivial component of \eqref{eq:violation_grad} is the Jacobian $\frac{\partial  \mathcal{O}_{\bm{p}}}{\partial \bm{y}}$. This information represents backpropagation through the lower-level problem, which can be obtained from one of the differentiable solvers discussed in the previous section, while automatic differentiation in PyTorch \citep{paszke2017automatic} is sufficient to complete the remaining chain rule calculations in \eqref{eq:violation_grad}.
We also define a function which reduces $ \| \nu(\bm{y}) \|^2$ by performing a gradient descent step of size $\gamma$:
\begin{equation}
% \vspace{-2pt}
    \label{eq:grad_step_fn}
       \cG(\bm{y}) \coloneqq \bm{y} - \gamma \nabla_{\bm{y}}\| \nu(\bm{y}) \|_2^2.
%\vspace{-6pt}
\end{equation}
Importantly, this function can also be rendered differentiable by leveraging functionality for \emph{back-propagating gradient calculations} (in this case, equation \eqref{eq:violation_grad}) in automatic differentiation libraries such as PyTorch. The result of function \eqref{eq:grad_step_fn} is generally infeasible to $\eqref{eq:ERM_c1}$, and this can be addressed by a (differentiable) projection back onto $\cC_{\bm{p}}$, completing one step of an end-to-end differentiable \emph{Coupling Constraint Correction} routine:
\vspace{-8pt}
\begin{equation}
    \label{CCC-pgd}
       \bm{y}_{k+1} = \bm{\Pi}_{\mathcal{C}_{\bm{p}}} \left( \bm{y}_k - \gamma \nabla \|\nu(\bm{y}_k)\|_2^2 \right).
\end{equation}
\vspace{-8pt}
Letting $\bm{y}_0 = \cN_{\theta}(\bm{p})$, our architecture for $\hat{\cB}_{\theta}$ corrects the neural network with $m$ steps of \eqref{CCC-pgd}. Explicitly, 
\begin{equation}
    \label{eq:architecture}
       \hat{\cB_{\theta}}(\bm{p}) =\left[ \left( \bm{\Pi}_{\mathcal{C}_{\bm{p}}} \circ \cG     \right)^m \circ \cN_{\theta} \right] (\bm{p}).
\end{equation}
By construction, this model is end-to-end differentiable, and maintains feasibility to \eqref{eq:ERM_c2} and  \eqref{eq:ERM_c1} while iterating toward satisfaction of \eqref{eq:ERM_c1}. Furthermore, the process \eqref{CCC-pgd} can be recognized as the classical projected gradient descent method on $\|\nu(\bm{y})\|_2^2$. This method is well-known to converge to local minima of convex and nonconvex functions, provided certain conditions on those functions' properties, their feasible set as well as the stepsize \citep{beck2017first}. The correction routine \eqref{CCC-pgd} is illustrated in Figure \ref{fig:illustration}, in which blue arrows represent gradient steps \eqref{eq:grad_step_fn}, and the alternating green arrows represent projections \eqref{C_projection} back onto $\cC_{\bm{p}}$. The entire chain of operations is composed with the neural network $\cN_{\theta}$, and unrolled in backpropagation to update its predictions.

%\subsection{Learning Bilevel Optimization with Coupling Constraints}
%So far, this section has proposed a method for learning to solve a parametric bilevel problem \eqref{UL} in the absence of its coupling constraint \eqref{UL-coupling}, in addition to a differentiable mechanism \eqref{CCC-pgd} to correct violations of \eqref{UL-coupling} incurred by a predicted solution. An overall method for learning to solve the full version of problem \eqref{UL} results from combining those two methods. It may be viewed as an adaptation of Algorithm \ref{alg:no_coupling}, in which the projection at line 4 is replaced with $m$ steps of projected gradient descent \eqref{CCC-pgd} on the coupling violation \eqref{eq:CCV}.  To enhance convergence to feasible solutions, a modified loss function is also adopted. 

\noindent \textbf{Penalty loss function.} The proposed architecture builds on the concept of a gradient-based constraint correction mechanism, popularized by the DC3 algorithm  \citep{donti2021dc3}.  Gradient descent methods are not guaranteed to converge in general, and it follows that DC3 corrections are not guaranteed to yield zero violations. However, as noted in \citep{donti2021dc3}, when \emph{initialized} close to an optimum, they are highly effective in
practice \citep{busseti2019solution, lee2019first}. Following \cite{donti2021dc3}, we use a penalty loss:
\begin{equation}
    \label{eq:soft_loss}
         \cL^{\textsc{soft}}_{\bm{p}_{(i)}}(\hat{\bm{y}},\hat{\bm{z}}) \coloneqq \cL_{\bm{p}_{(i)}}(\hat{\bm{y}},\hat{\bm{z}}) + \lambda \|\nu(\hat{\bm{y}})\|_2^2 .
\end{equation}
This leads to upper-level predictions from $\cN_{\theta}$ which are initialized close to satisfying the coupling constraint~\eqref{UL-coupling}, so that its end-to-end training with the differentiable correction~\eqref{CCC-pgd} tends to produce feasible solutions to the full problem \eqref{UL}.

\noindent \textbf{Training Routine.}
An overall training scheme is summarized in Algorithm \ref{alg:NCC}, in terms of one epoch of stochastic gradient descent. Each data input sample $\bm{p}_{(i)}$ represents a distinct instance of problem \eqref{UL}, for which an initial solution estimate $\hat{\bm{y}}$ is predicted at line 3. The sequence \eqref{CCC-pgd} makes up lines $4$-$8$. For each of $m$ \emph{correction steps}, $\hat{\bm{y}}$ is iteratively refined by taking a step toward feasibility of the coupling constraint \eqref{UL-coupling}. The gradient $\hat{\bm{g}}$ at line 6 is calculated per equation \eqref{eq:violation_grad}. The loss \eqref{eq:soft_loss} is then evaluated w.r.t. the refined estimate $\hat{\bm{y}}$ and its lower-level pair.  If needed, at test time, more than $m$ iterations may be applied in the correction routine.
Note that Algorithm \ref{alg:NCC} describes only the \emph{forward pass} of the training routine. Line 11 encapsulates backpropagation through all components of the model, implemented with a combination of automatic differentiation and implicit differentiation via the differentiable solvers, as described above.
%%\begin{wrapfigure}{R}%%{0.48\textwidth}
%%%\vspace{-10pt} % Adjust as needed
%%\begin{minipage}{0.46\textwidth}
%%\setlength{\intextsep}{0pt} % %%Reduce space below the wrapfigure
%%\setlength{\columnsep}{10pt} % %%Optional: adjust spacing between %%text and figure

\begin{wrapfigure}{R}{0.5\textwidth} % Adjusted width for half-page
\vspace{-20pt} % fine-tune vertical spacing
\begin{minipage}{0.48\textwidth} % ensures proper sizing within wrapfigure
\begin{algorithm}[H]
   \caption{Learning Bilevel Optimization \\ with Coupling}
   \label{alg:NCC}
\begin{algorithmic}[1]
   \STATE {\bfseries Input:} parameters $\{\bm{p}_{(i)}\}_{i=1}^N$, \\weights $\bm{\theta}$, \\learning rate $\alpha$, correction stepsize $\gamma$  
   \FOR{$i = 1$ {\bfseries to} $N$}
   \STATE  $\hat{\bm{y}} \gets \cN_{\bm{\theta}}(\bm{p}_{(i)})$ 
   \FOR{$k=1$ {\bfseries to} $m$}
   \STATE $\hat{\bm{g}} \gets \nabla \left\| \textsc{ReLU}\left( \bm{U}_{\bm{p}_{(i)}}( \hat{\bm{y}},\cO_{\bm{p}_{(i)}}(\hat{\bm{y}}) ) \right) \right\|^2$
   \STATE $\hat{\bm{y}} \gets \hat{\bm{y}} - \gamma \cdot \hat{\bm{g}}$
   \STATE $\hat{\bm{y}} \gets {\mathbf{\Pi}}_{\mathcal{C}_{\bm{p}_{(i)}}}(\hat{\bm{y}})$
   \ENDFOR
   \STATE $\hat{\bm{z}} \gets  \cO_{\bm{p}_{(i)}}(\hat{\bm{y}})$
   \STATE $\bm{g} \gets \nabla_{\theta} \cL^{\textsc{soft}}_{\bm{p}_{(i)}}(\hat{\bm{y}},\hat{\bm{z}})$
   \STATE $\bm{\theta} \gets \bm{\theta} - \alpha \cdot \bm{g}$
   \ENDFOR
\end{algorithmic}
\end{algorithm}
\end{minipage}
\vspace{-6pt}
\end{wrapfigure}
%%\end{minipage}
%%\vspace{-10pt} 
%%\end{wrapfigure}

%\paragraph{Distinctions from DC3.} The idea of a gradient-based constraint correction in Learning to Optimize was first popularized by the DC3 method \cite{donti2021dc3}, aimed at nonlinear optimization problems of a single level. This paper generalizes that approach, as part of an overall framework for learning \emph{bilevel} optimization. \rev{put this earlier in the paper and expand on the rationale: partitioning constraints into coupling and non-coupling; recognizing non-coupling as the easier ones, satisfy them with projections and iterate toward satisfying the other. Diff opt is an essential component not needed in the og DC3.}

\iffalse
\begin{algorithm}[tb]
   \caption{Coupling Constraint Correction}
   \label{alg:CCC}
\begin{algorithmic}[1]
   \STATE {\bfseries Input:} parameters $\bm{p}$, predicted variables $\hat{\bm{y}}$, stepsize $\gamma$, $\#$ steps $m$  
   \FOR{$i=1$ {\bfseries to} $m$}
   \STATE $\hat{\bm{y}} \gets {\mathbf{\Pi}}_{\mathcal{C}_{\bm{p}}}(\hat{\bm{y}})$  \label{line:CCC_project} \\
   \STATE $\hat{\bm{z}} \gets  \cO_{\bm{p}}(\hat{\bm{y}})$  \label{line:CCC_opt}\\
   \STATE $\bm{g} \gets \nabla \| \left[ \bm{U}_{\bm{p}}( \hat{\bm{y}},\hat{\bm{z}} ) \right]_+ \|^2 \;\;$ \label{line:CCC_grad}  \\
   \STATE $\hat{\bm{y}} \gets \hat{\bm{y}} + \gamma \cdot \bm{g}$ \label{line:CCC_step} 
   \ENDFOR
\end{algorithmic}
\end{algorithm}
\fi

\subsection{Satisfying Lower-Level Constraints}
\label{sec:LowerLevelConstraints}
So far, this section has assumed that the lower-level problem \eqref{LL} had at least one feasible solution for any $\bm{p}$ and $\bm{y} \in \cC_{\bm{p}}$ (see the comments on existence and uniqueness in Section \ref{sec:Setting}). Otherwise, the optimization problem at line $5$ of Algorithm \ref{alg:NCC} may be infeasible. In such cases, we address this issue by reformulating the overall bilevel problem so that the property is satisfied. To do so, we simply identify the lower-level constraints that prevent feasibility and \emph{lift} them to the upper level where they become coupling constraints.
Let the lower level's feasible set be partitioned as $\cS_{\bm{p}}(\bm{y}) = \cF_{\bm{p}}(\bm{y}) \cup \cF^{C}_{\bm{p}}(\bm{y})$, such that $\cF_{\bm{p}}(\bm{y})$ is nonempty for all $\bm{p}$ and $\bm{y} \in \cC_{\bm{p}}$. We reformulate problems \eqref{UL} and \eqref{LL} so that $\cS_{\bm{p}}(\bm{y})$ is replaced by $\cF_{\bm{p}}(\bm{y})$ in problem \eqref{LL}, while $\bm{z} \in \cF^{C}_{\bm{p}}(\bm{y})$ is promoted to the upper level problem \eqref{UL}. As an upper-level constraint relating $\bm{y}$ and $\bm{z}$, it is absorbed into the \emph{coupling constraints} \eqref{UL-coupling}. This technique is applied in the experiments of Section \ref{sec:ControlCoDesign}, and detailed in Appendices \ref{app:TT} and \ref{app:HVAC}.

\section{Experiments}
\label{sec:Experiments}

In this section, we evaluate the proposed methods' effectiveness in learning to solve several parametric bilevel optimization problems. In order to measure the optimality of its learned solutions, we require true optimal solutions for comparison. In light of Section \ref{sec:RelatedWork}, finding those solutions can pose significant challenges in itself. Most known methods either rely on exploiting problem-specific structure or otherwise employ metaheuristic methods which lack guarantees. Even if optimal solutions are found by such methods, their optimality often cannot be certified \citep{sinha2017review}. For these reasons, Section \ref{sec:bilevel_QP} begins by learning solutions to small-scale synthetic problems for which open-source solvers can provide \emph{certified optimal solutions for comparison}. Then, Section \ref{sec:ControlCoDesign} extends the evaluation to more complex bilevel programs in engineering design. Those problems are significantly more difficult due to their much larger size and more complex forms, which include \emph{nonconvex optimization} at the lower level. In those cases, we evaluate our learned solutions against those of a Particle Swarm Optimization (PSO), a metaheuristic framework commonly used in design optimization \citep{asaah2021optimal, hou2015optimized}.

\noindent \textbf{Evaluation criteria and conventions.}
We evaluate the ability of Algorithm \ref{alg:NCC} to perform the training task specified in Equation \ref{eq:ERM}. Recall that constraints \eqref{eq:ERM_c1} and \eqref{eq:ERM_c2} are ensured by the construction of Algorithm \ref{alg:NCC}; thus, our two main evaluation criteria are the \emph{objective values} \eqref{eq:ERM_objective}, and potential \emph{coupling constraint violation} \eqref{eq:CCV}. When true optimal solutions $({\bm{y}}^{\star}, {\bm{z}}^{\star})$ are available, we report the \emph{relative optimality gap} relative to learned solutions $({\hat{\bm{y}}}, {\hat{\bm{z}}})$, which we define $\left| \nicefrac{ \mathcal{L}(\hat{\bm{y}},\hat{\bm{z}}) - \mathcal{L}(\bm{y^{\star},z^{\star}}) }{\mathcal{L}(\bm{y^{\star},z^{\star}})} \right|$
and illustrate \emph{in blue} throughout. When the true optima are not known, we instead report the nominal objective value of the learned solutions, \emph{in green}, alongside the solution produced by a baseline method for comparison. All metrics are reported on average over the respective test set. When a metric should ideally converge to zero, its standard deviation is also reported. Additional implementation details for each experiment, can be found in Appendices \ref{app:BQP}, \ref{app:TT}, \ref{app:HVAC} .

\begin{table*}[!t]
\centering
\begin{minipage}[t]{0.42\textwidth}
\centering
\resizebox{\linewidth}{!}{
\begin{tabular}{l | r r r | r}
    \toprule
        \textbf{BQP}     & 
        $\left| \frac{ \mathcal{L}(\hat{\bm{y}},\hat{\bm{z}}) - \mathcal{L}(\bm{y^{\star},z^{\star}}) }{\mathcal{L}(\bm{y^{\star},z^{\star}})} \right|$ &
        $\| \nu(\bm{y}) \|_2$  &
        Time (s) &
        YALMIP (s) \\
    \midrule
     $3 \times 2$  & $9.2e-4 \pm 2e-3$ & $5.9e-3  \pm 1e-2$  & $\bm{6.4e-2}$ & $0.12$\\ 
     $6 \times 4$  & $2.0e-3 \pm 5e-3$ & $2.8e-4 \pm 9e-4$ & $\bm{6.4e-2}$ & $1.2$\\  
     $9 \times 6$  & $1.1e-2 \pm 1e-2$ & $4.0e-5 \pm 7e-4$ & $\bm{6.7e-2}$ & $10.3$\\
    \bottomrule
\end{tabular} 
}
\caption{BQP problems, Test Set Average}
\label{table:BQP_metrics}
\end{minipage}
\vspace{-8pt}
\hfill
\begin{minipage}[t]{0.57\textwidth}
\centering
\resizebox{\linewidth}{!}{
\begin{tabular}{l | rr | rr | rr}
    \toprule
        \multicolumn{1}{c}{\textbf{Model}}     & 
        \multicolumn{2}{c}{$\cL(\bm{y})$} &
        \multicolumn{2}{c}{$\| \nu(\bm{y}) \|_2$}  &
        \multicolumn{2}{c}{Time (s)}\\
    \cmidrule{2-3}\cmidrule{4-5}\cmidrule{6-7}
    & {TT} & {HVAC}
    & {TT} & {HVAC}
    & {TT} & {HVAC}\\
    \midrule
     Learned (Alg. \ref{alg:NCC})  & $\bm{0.122}$ & $\bm{1.23}$ & $\bm{1.4e-2} \pm 2e-2$ & $3.0e{-2} \pm 5e-2$ & $\bm{2.7e-2}$  & $\bm{7.70e-2}$\\ 
     PSO (baseline) & $0.140$  & $1.46$ & $1.6e-2 \pm 2e-2$ & $3.0e{-2} \pm 2e-2$ & $1268.7$ & $1055.9$\\  
    \bottomrule
\end{tabular} 
}
\caption{On Control Co-Design, Test Set Average}
\label{table:codesign_metrics}
\end{minipage}
\vspace{-8pt}
\end{table*}

\vspace{-4pt}
\subsection{Preliminary Problem: Learning Bilevel Quadratic Programming}
\label{sec:bilevel_QP}
We begin the experimental evaluation on a relatively simpler class of bilevel problems, the Bilevel Quadratic Programs (BQP). Both their upper and lower-level problems contain convex quadratic objective functions and only linear constraints. The upper and lower-level problems are respectively:
\begin{minipage}{0.54\textwidth}
\vspace{-8pt}
\begin{subequations}
    \label{QP-UL}
    \begin{align}
         \cB(\bm{c},\bm{d}) = \underset{\bm{y}}{\argmin} \;\;\;\; & \frac{1}{2}\bm{y}^T \bm{Q} \bm{y} + \bm{c}^T\bm{y} + \bm{d}^T\bm{z} + q  \\
         \textit{s.t.}\;\;\;\; & \bm{A}\bm{y} \leq \bm{b} + \bm{E} \bm{z} \\
         & \bm{z} \in \cO(\bm{y}) ,
    \end{align}
\end{subequations}
\vspace{-16pt}
\end{minipage}
\hfill
\begin{minipage}{0.45\textwidth}
\vspace{-8pt}
\begin{subequations}
    \label{QP-LL}
    \begin{align}
         \cO(\bm{y}) = \underset{\bm{z}}{\argmin} \;\; & \frac{1}{2}\bm{z}^T \bm{H} \bm{z} + \bm{e}^T\bm{z} + \bm{f}^T\bm{y} + g  \\
         \textit{s.t.}\;\;\;\; & \bm{F}\bm{z} \leq \bm{h} + \bm{G} \bm{y} .
    \end{align}
\end{subequations}
\vspace{-16pt}
\end{minipage}

\iffalse
\begin{subequations}
    \label{QP-UL}
    \begin{align}
         \cB(\bm{c},\bm{d}) = \underset{\bm{y}}{\argmin} \;\;\;\; & \frac{1}{2}\bm{y}^T \bm{Q} \bm{y} + \bm{c}^T\bm{y} + \bm{d}^T\bm{z} + q  \\
         \textit{s.t.}\;\;\;\; & \bm{A}\bm{y} \leq \bm{b} + \bm{E} \bm{z} \\
         & \bm{z} \in \cO(\bm{y}) ,
    \end{align}
\end{subequations}
where
\begin{subequations}
    \label{QP-LL}
    \begin{align}
         \cO(\bm{y}) = \underset{\bm{z}}{\argmin} \;\; & \frac{1}{2}\bm{z}^T \bm{H} \bm{z} + \bm{e}^T\bm{z} + \bm{f}^T\bm{y} + g  \\
         \textit{s.t.}\;\;\;\; & \bm{F}\bm{z} \leq \bm{h} + \bm{G} \bm{y} .
    \end{align}
\end{subequations}
\fi
We train a neural network to learn its solutions as a function of the upper-level linear objective coefficients. That is, in the notation of Section \ref{sec:Setting}, we have $\bm{p} \coloneqq (\bm{c},\bm{d})$.

\noindent \textbf{Experimental details.}
When problems \eqref{QP-UL} are sufficiently small, they can be solved using a mixed-integer programming reformulation. In this experiment, \emph{we consider problems within a range of sizes such that ground-truth optimal solutions to the test-set problems can be computed within a reasonable time}. \emph{This provides an initial setting} in which certified optimal solutions can be obtained, against which the optimality of our learned solutions can be measured objectively. We use the open-source YALMIP package to solve the instances by replacing their lower-level problem with KKT conditions, and solving the resulting single-level MIP with a branch-and-bound method \citep{lofberg2004yalmip}.

We refer to a BQP problem with $m$ and $n$ variables at the upper and lower levels, respectively, as having size $m \times n$. Three sets of BQP problems having size $3 \times 2$, $6 \times 4$, and $9 \times 6$ are randomly generated - first by drawing the elements of each matrix $\bm{A},\bm{E},\bm{F},\bm{G}$ and each vector $\bm{b}, \bm{e}, \bm{f}, \bm{h}$ from a uniform distribution $U(0,1)$. Positive-semidefinite $\bm{Q}$ and $\bm{H}$ are constructed by self-multiplication of such a matrix. Individual problem instances are generated by drawing vectors of linear objective coefficients $(\bm{c}, \bm{d})$ also from $U(0,1)$. Each set is divided into validation and test portions, numbering $1000$ each. The prediction model $\hat{\cB}_{\theta}$ is a $5$-layer network followed by $20$ steps of Algorithm \ref{alg:NCC}. Solution and differentiation of problem \eqref{QP-LL} are implemented in  \emph{cvxpylayers} \citep{agrawal2019differentiable}.

\begin{figure*}
\centering
\includegraphics[width=0.9\textwidth]{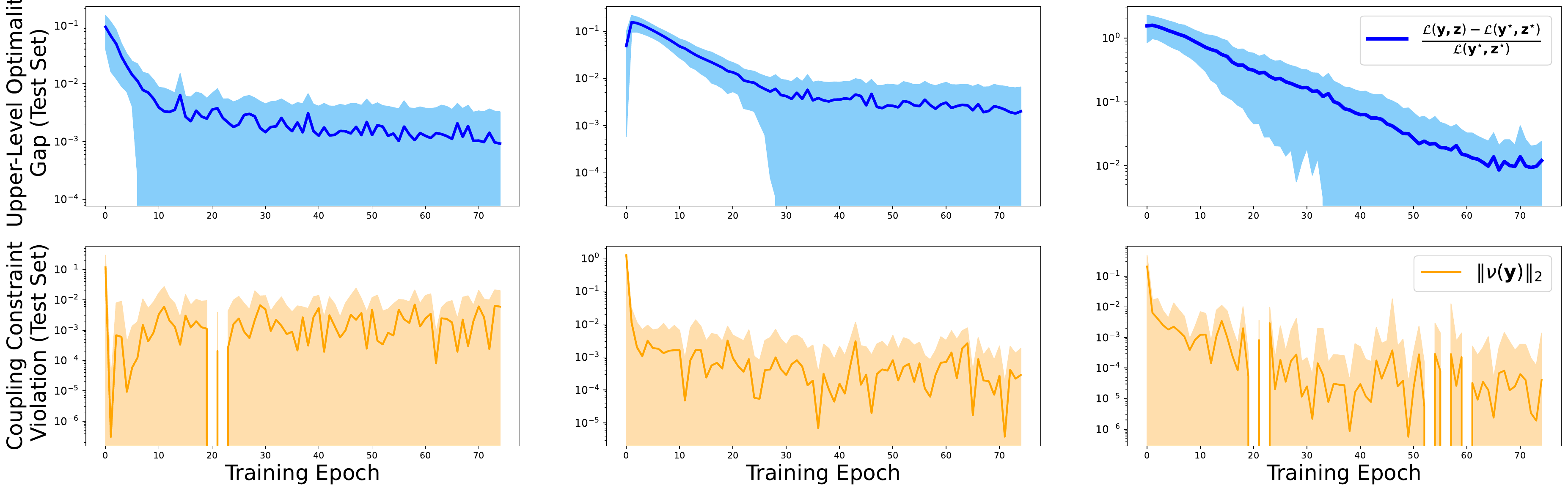}
\caption{Optimality gap and coupling violation, on different-sized BQP problems from left to right: $3 \times 2$, $6 \times 4$, $9 \times 6$. Shows mean and standard deviation over the test set, at each training epoch. }
\label{fig:QP_curves}
\vspace{-12pt}
\end{figure*}

\noindent \textbf{Results.}
Figure \ref{fig:QP_curves} illustrates the two main evaluation metrics throughout Algorithm \ref{alg:NCC}, in terms of mean and standard deviation over the test set for each parametric BQP. In each case, the relative optimality gap   (in blue)   is reduced by $2$ orders of magnitude over $75$ epochs to a value between $10^{-3}$ and $10^{-2}$. The coupling violation \eqref{eq:CCV} is rapidly diminished in the first epoch and then generally bounded below $10^{-2}$ within a full standard deviation throughout training. Test set metrics are also reported in Table \ref{table:BQP_metrics}. Together, these results demonstrate the ability of Algorithm \ref{alg:NCC} to learn bilevel optimization with negligible error on small-scale BQP problems. Beyond the problem sizes considered here, starting with $12 \times 9$, the time taken by YALMIP to fully solve the test instances \emph{becomes intractable, thus we cannot benchmark against certified optimal solutions on larger instances}.

\vspace{-6pt}
\subsection{Learning Optimal Control Co-Design}
\label{sec:ControlCoDesign}

This section introduces more challenging bilevel optimizations, in terms of both their size and their form. Optimal Control Co-Design is a bilevel problem setting in which an engineering system is designed to optimize an economic objective at the upper level, subject to conditions on its behavior under a known optimal control policy, which forms its lower level. The problems described below cannot hope to be solved by conventional methods with certificates of optimality; as an alternative, we compare our learned solutions to the results of a PSO-based metaheuristic method. The PSO framework is commonly applied to problems that lack efficient solution methods, and thus is a favored tool in design optimization \citep{asaah2021optimal,hou2015optimized}. Details of the PSO baseline methods, along with illustrations of their results, can be found in Appendix \ref{app:PSO_details} . %Throughout the section, variables $\bm{u}$ and $\bm{x}$ quantify the control actions and resulting states of a system.

\vspace{-6pt}
\subsubsection{Nonconvex Bilevel Optimization: Control Co-design of a Nonlinear  System}
\label{sec:twotank}

We consider a nonlinear control problem, in which two connected tanks are controlled by a single pump and a two-way valve. The system is a simplified model of a pumped-storage hydroelectricity, which is a form of  energy storage used by electric power systems for load balancing. The system dynamics are described by nonlinear ODE's $\dot{\bm{x}} = f(\bm{x}^{(i)}, \bm{u}^{(i)}, \bm{y})$:\\[-2pt]
\begin{minipage}{0.48\textwidth}
\vspace{-0.5\baselineskip}
\begin{subequations}
\label{eq:twotank_ODE_1}
\begin{align}
    \dot{x_1} = y_1 (1-u_2) u_1 - y_2 \sqrt{x_1},
\end{align}
\end{subequations}   
\vspace{-\baselineskip}
\end{minipage}
\hfill
\begin{minipage}{0.48\textwidth}
\vspace{-0.5\baselineskip}
\begin{subequations}
\label{eq:twotank_ODE_2}
\begin{align}
    \dot{x_2} = y_1 u_2 u_1 + y_2 \sqrt{x_1} - y_2 \sqrt{x_2},
\end{align}
\end{subequations}
\vspace{-\baselineskip}
\end{minipage}

\noindent in which $x_1$, $x_2$ are the water levels in each tank. Control actions consist of $u_1$ and $u_2$, which are the pump modulation and valve opening. The nonlinear optimal control problem \eqref{eq:twotank_LL_discrete} seeks the control policy which minimizes energy expended to reach a desired terminal state $\bm{p}$. The function \texttt{ODESolve} represents Euler discretization of (\ref{eq:twotank_ODE_1}, \ref{eq:twotank_ODE_2}) over $N=20$ frames to a final time $T$. This yields new variables $\bm{x}=[\bm{x}^{(1)}, \ldots, \bm{x}^{(N)}] \in \mathbb{R}^{N \times 2}$, $\bm{u} = [\bm{u}^{(1)}, \ldots, \bm{u}^{(N)}] \in \mathbb{R}^{N \times 2}$ bound by a sequence of nonlinear (i.e., \emph{nonconvex}) equality constraints \eqref{eq:twotank_LL_discrete_controls}  for $1 \leq i < N$, while  $dt = \frac{T}{N}$.
\begin{figure}[h!]
\centering
\begin{minipage}{0.49\textwidth}
\vspace{-1.5\baselineskip}
\begin{subequations}
\label{eq:twotank_UL}
\begin{align}
    \mathcal{B}({\bm{p}}) = {\argmin}_{\bm{y}}\;\;\; & \bm{v}^T \bm{y} \label{eq:twotank_UL_obj} \\
    \textit{s.t.}\;\;\; 
    & \bm{x},\bm{u} = \mathcal{O}_{\bm{p}}(\bm{y})  \label{eq:twotank_UL_opt} \\ 
    & \bm{y}_{min} \leq \bm{y} \leq \bm{y}_{max} \label{eq:twotank_bounds} \\
    & \bm{x}^{(N)} = \bm{p},  \label{eq:twotank_UL_endpt}
\end{align}
\end{subequations}
\vspace{-2\baselineskip}
\end{minipage}
\hfill
\begin{minipage}{0.49\textwidth}
\vspace{-1.5\baselineskip}
\begin{subequations}
\label{eq:twotank_LL_discrete}
\begin{align}
    \mathcal{O}_{\bm{p}}(\bm{y}) 
    & \coloneqq \underset{0 \leq \bm{x}, \bm{u} \leq 1}{\argmin}   
    \;\;\; \sum_{k=1}^N \| \bm{u}^{(k)} \|_2^2 \\
    \textit{s.t.}\;\;\;
    & \bm{x}^{(N)} = \bm{p} \label{eq:twotank_LL_discrete_ref}\\
    & \bm{x}^{(i+1)} = \text{ODESolve}\big(f(\bm{x}^{(i)}, \bm{u}^{(i)}, \bm{y})\big). 
    \label{eq:twotank_LL_discrete_controls}
\end{align}
\end{subequations}
\vspace{-2\baselineskip}
\end{minipage}
\end{figure}
The upper-level problem \eqref{eq:twotank_UL} seeks to optimize the design of such a system in terms of its overall cost $\mathcal{L}(\bm{y},\bm{z}) \coloneqq \bm{v}^T \bm{y}$, treating the inlet and outlet valve coefficients $\bm{y} = [y_1, y_2]$ as free design parameters. A feasible design demands that upper and lower bounds on each element of $\bm{y}$ are satisfied per \eqref{eq:twotank_bounds}. Additionally, the parametric end-state constraint \eqref{eq:twotank_LL_discrete_ref} is duplicated at the upper level in \eqref{eq:twotank_UL_endpt} to emphasize its coupling role. The full reformulation per Section \ref{sec:LowerLevelConstraints} is detailed in Appendix \ref{app:TT}.  The initial condition is $\bm{x}^{0}=\bm{0}$. Taken together, the coupled problems (\ref{eq:twotank_UL}, \ref{eq:twotank_LL_discrete}) seek the parameters $\bm{y}$ which yield the minimal-cost system design that can be controlled to state $\bm{p}$ by time $T$ under its control policy \eqref{eq:twotank_LL_discrete}. In this experiment, our model $\hat{\cB}_{\theta}$ is trained to perform a fast and accurate approximation to this design problem for any such scenario specified by a given $\bm{p} \in [0,1]$. In the notation of Section \ref{sec:Setting},  $\bm{z} \coloneqq (\bm{x},\bm{u})$.

\noindent \textbf{Experimental details.}
We consider an experiment in which $T=10s$ and $N=20$, $\bm{c}_{min}=0$, and $\bm{c}_{max}=\frac{1}{3}$. Problem instances correspond to reference states $\{\bm{p}\}$, which are randomly generated from $U(0,1)$, but with $p_1 < p_2$ to ensure feasibility of problem \eqref{eq:twotank_UL}. They are partitioned into training, validation, and test sets of sizes $10000$, $1000$, $1000$.  The model $\hat{\cB}_{\theta}$ consists of an $8$-layer network ${\cN}_{\theta}$ followed by Algorithm \ref{alg:NCC}.  To implement Algorithm \ref{alg:NCC} is nontrivial, as it requires a differentiable solution of the nonconvex lower-level programs \eqref{eq:twotank_LL_discrete}. For this, we employ the differentiable model predictive control solver of  \cite{amos2018differentiable}, which differentiates problem \eqref{eq:twotank_LL_discrete} implicitly via the KKT conditions of the final convex subproblem of a sequential quadratic programming.

\noindent \textbf{Results.}
Figure \ref{fig:tt_training_curve} (top two) illustrates the value of the design objective $\mathcal{L}( \mathbf{y}, \mathbf{z} )$, as well as the coupling violation \eqref{eq:CCV}, over the test set throughout training.  Overall metrics are also found in Table \ref{table:codesign_metrics} under the header TT. Despite investing an average of 1286.7s of solution time per instance, the PSO baseline produces design solutions with $15$ percent higher cost than those learned by Algorithm \ref{alg:NCC}, which infers solutions in $0.027$s on average. At the same time, it attains nearly identical satisfaction of the coupling constraint, on average. \emph{This result is significant, since it demonstrates an ability to learn nonconvex bilevel optimization with high accuracy}. 
\iffalse
\begin{figure*}
\centering
\includegraphics[width=0.43\textwidth]{images/tt_curves.pdf}
\includegraphics[width=0.41\textwidth]{images/HVAC_curves.pdf}
% \includegraphics[width=0.9\columnwidth]{images/portfolio_objective.pdf}
\caption{Test set metrics, learning control co-design of the two-tank system (left) and HVAC system (right).  }
\label{fig:tt_training_curve}
\end{figure*}
\fi

\begin{figure}
\centering
\includegraphics[width=0.43\textwidth]{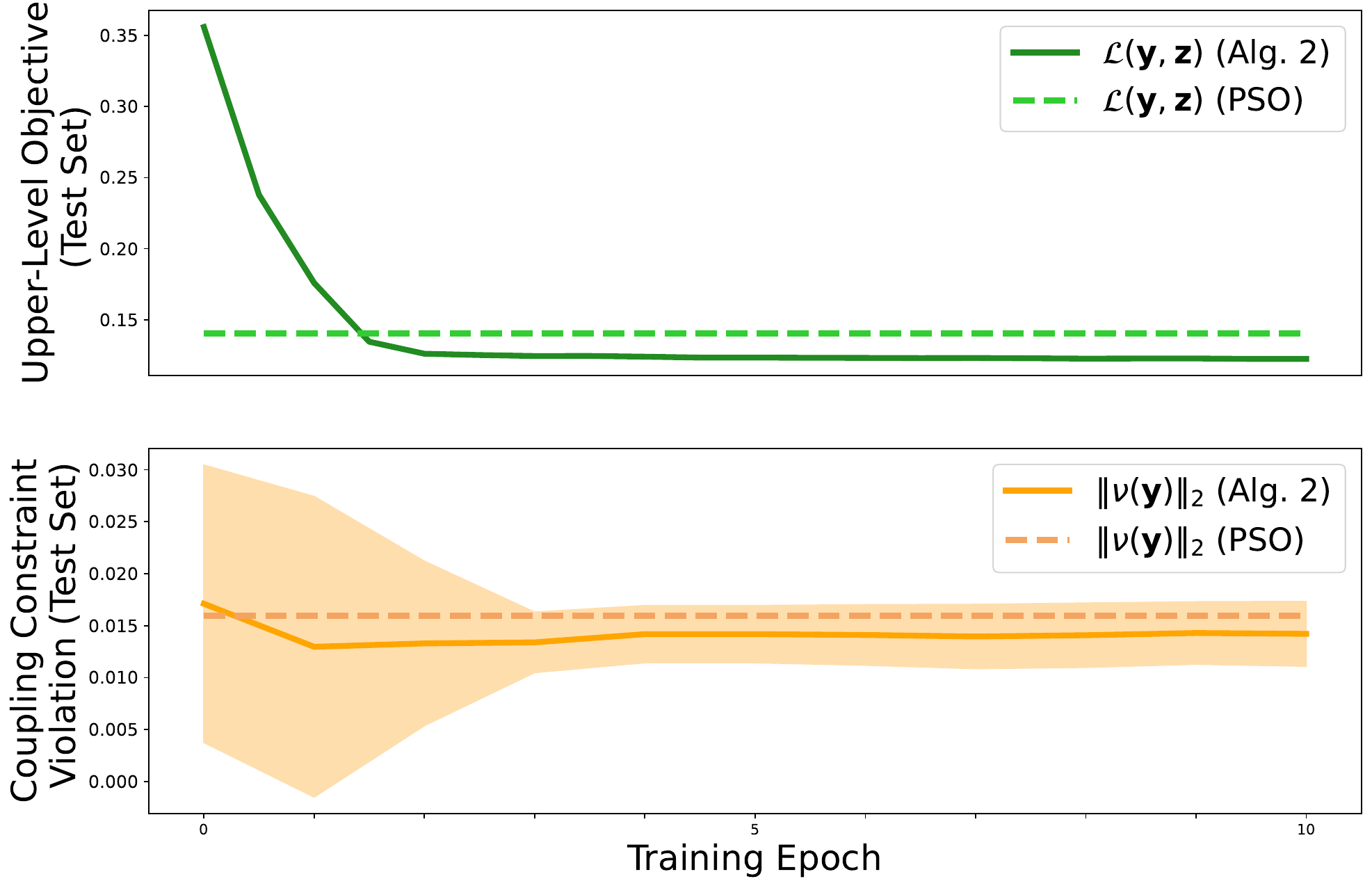}
\vspace{-6pt}
\includegraphics[width=0.43\textwidth]{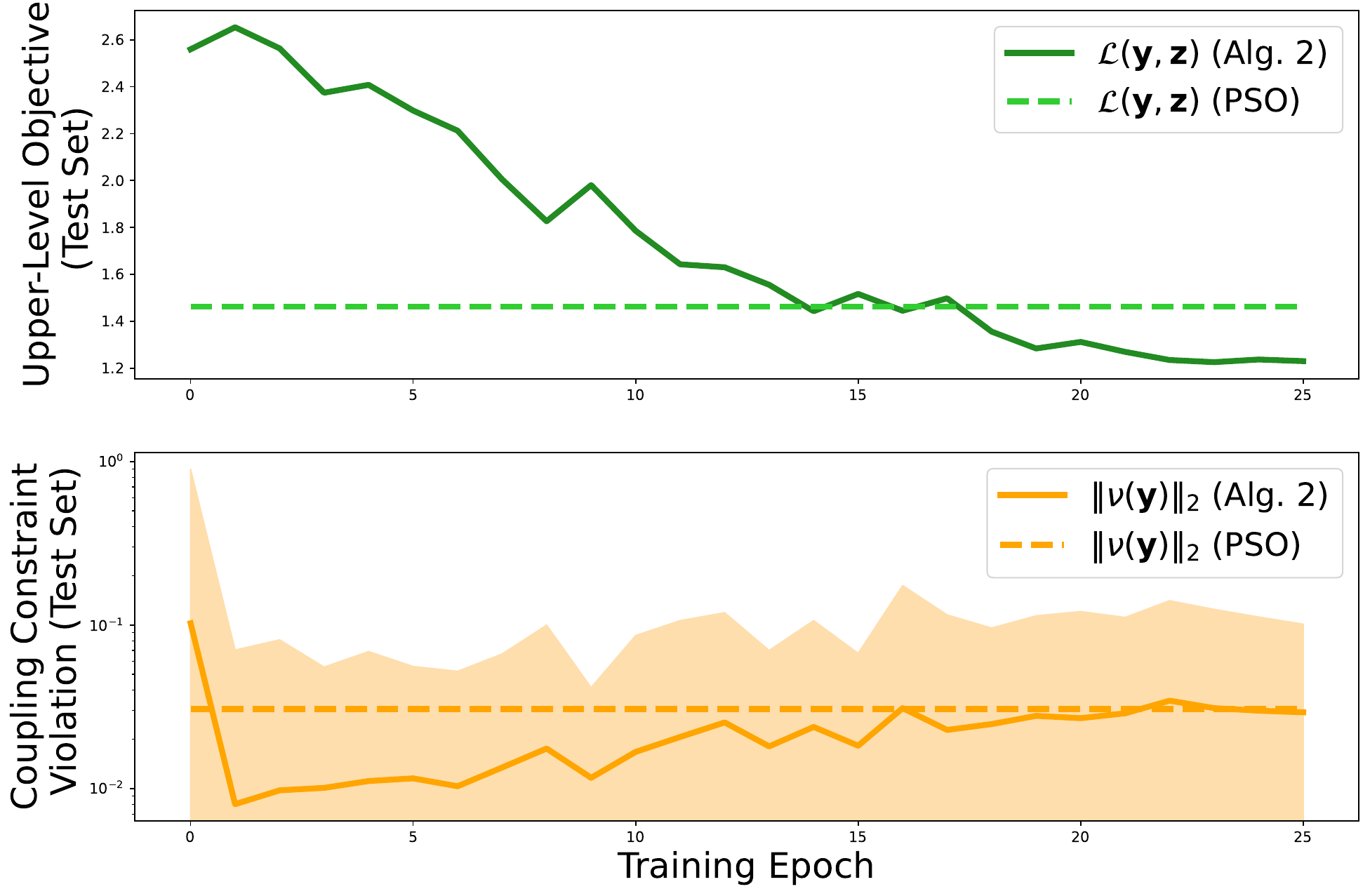}
\caption{ Test set metrics, learning control co-design of the two-tank (left) / HVAC (right) systems. }
\label{fig:tt_training_curve}
\vspace{-12pt}
\end{figure}
\subsubsection{Control Co-Design of a Building HVAC System}
\label{sec:building}
Finally, we consider the design and control of the heating, ventilation, and air conditioning (HVAC) system in a building. The control problem minimizes energy use while maintaining indoor temperature within prescribed bounds~\eqref{eq:building}. The building consists of $2$ zones, thermally connected to each other and the outside environment by a matrix of conductivity coefficients $\bm{A}$. State variables $\bm{x}=[\bm{x}^{(1)}, \ldots, \bm{x}^{(N)}] \in \mathbb{R}^{N \times 8}$ consist of the temperatures of each zone's floor, walls, indoor air and exterior facade at each timestep $k$. Control actions $\bm{u} = [\bm{u}^{(1)}, \ldots, \bm{u}^{(N)}] \in \mathbb{R}^{N \times 2}$ induce heat flows into each zone, which affect the temperature states $\bm{x}$ via the actuator design variables $\bm{Y} \in \mathbb{R}^{8 \times 2}$. The states are also affected by random disturbances $\bm{d}$, which include heat transfer from occupants and solar irradiation. Thermal constraints \eqref{eq:HVAC_LL_bds} demand that the indoor air temperature must remain within prescribed time-varying bounds $(\underline{\bm{p}}, \overline{\bm{p}})$. This condition couples the upper and lower problems, which we emphasize by its duplication at \eqref{eq:HVAC_UL_bds}. The design task \eqref{eq:building_UL} asks to learn $\bm{Y}$ which minimizes a linear cost function $\textit{Tr}(\bm{V}^T \bm{Y})$ while allowing the system to be maintained within the thermal bounds using optimal control. In the notation of Section \ref{sec:Setting}, we have $\bm{p} \coloneqq (\underline{\bm{p}}, \overline{\bm{p}})$, $\bm{y} \coloneqq \bm{Y}$ and $\bm{z} \coloneqq (\bm{x},\bm{u},\bm{w})$:
%\begin{minipage}[t]{0.4\textwidth}
%\centering
\begin{figure}[h!]
\centering
\begin{minipage}{0.42\textwidth}
%\hspace*{-3.0em}
\vspace{-1.5\baselineskip}
\begin{subequations}
    \label{eq:building_UL}
    \begin{align}
         \cB(\bm{p}) = {\argmin}  & \textit{Tr}(\bm{V}^T \bm{Y})   \\
         \textit{s.t.}\;\; & \bm{x},\bm{u},\bm{w} = \mathcal{O}_{\bm{p}}(\bm{Y}) \\
          & \bm{Y} \geq \bm{0} \\ & \underline{\bm{p}}^{(k)} \leq \bm{w}^{(k)} \leq \overline{\bm{p}}^{(k)}  \label{eq:HVAC_UL_bds} 
    \end{align}
\end{subequations}
\vspace{-\baselineskip}
\end{minipage}
%\hspace*{1.0em}
\begin{minipage}{0.56\textwidth}
\vspace{-2\baselineskip}
\begin{subequations}
    \label{eq:building}
    \begin{align}
         \cO_{\bm{p}}(\bm{Y}) = \underset{\bm{x}, 0 \leq \bm{u} \leq 1, \bm{w}}{\argmin}  & \sum_{k \in \{1 \ldots N \}} \| \bm{u}^{(k)} \|_2^2   \\
         \textit{s.t.}\;\;  
          & \bm{w}^{(k)} = \bm{C}\bm{x}^{(k)} \\
          &  \underline{\bm{p}}^{(k)} \leq \bm{w}^{(k)} \leq \overline{\bm{p}}^{(k)}  \label{eq:HVAC_LL_bds} \\
            &\bm{x}^{(k+1)} = \bm{A}\bm{x}^{(k)} + \bm{Y}\bm{u}^{(k)} + \bm{E}\bm{d}^{(k)} 
            \end{align}
\end{subequations}
\vspace{-3\baselineskip}
\end{minipage}
\end{figure}
%\label{eq:LL}
%\vspace{-12pt}
%\end{minipage}

\noindent \textbf{Experimental details.}
Our experiment assumes $N=30$ steps, the problem instances with for the thermal bounds $\underline{\bm{p}}$ are generated from a $\beta$-random walk along with $ \overline{\bm{p}} = \underline{\bm{p}} + 2.0$, validation and test sets of size $10000$, $1000$, $1000$. A fixed disturbance pattern $\bm{d}$ is generated from the building control test suite in \texttt{NEUROMANCER} \citep{Neuromancer2023}. The model $\cB_{\theta}$ is a $6$-layer ReLU network followed by $10$ steps of Algorithm \eqref{alg:NCC}. Differentiable solution of     \eqref{eq:building} is implemented using \texttt{cvxpylayers}.

\noindent \textbf{Results.}
Figure \ref{fig:tt_training_curve} (right) illustrates the value of the design objective $\mathcal{L}( \mathbf{y}, \mathbf{z} )$, as well as the coupling violation \eqref{eq:CCV}, over the test set throughout training. Overall metrics are also found in Table \ref{table:codesign_metrics} under HVAC. While our learned designs incur nearly identical coupling constraint violations on average, they are achieved at about $19$ percent lower cost, and with orders of magnitude lower solving time. 

\vspace{-7pt}
\section{Conclusion and Limitations}
\label{sec:conclusion}
This paper has shown how the modern toolkit of differentiable optimization can be used to train machine learning models as fast and accurate approximators of parametric bilevel optimization with coupling constraints. Experiments on control co-design problems show that the proposed learning to optimize framework can accurately approximate even nonconvex programs. Multiple interesting avenues remain to extend the work. For example, extension to integer-valued decision variables would broaden its applicability. Second, when lower‑level objectives are nondifferentiable with respect to the upper variables, extensions based on surrogate gradients may be investigated. Addressing these points would further broaden the scope of this work.

% Future work may focus on extending the framework beyond the limitation which requires unique solutions at the lower level. Cases where the lower-level problem can not be differentiated with respect to upper-level variables should also be studied -  although differentiable surrogates for such problems have been developed, such settings are beyond the present scope. It is hoped that this work will help lead to increased interest in those topics, and in machine learning and differentiable optimization as part of an overall toolkit for tackling challenging bilevel optimization problems. 

\begin{ack}
% Use unnumbered first level headings for the acknowledgments. All acknowledgments

% Do {\bf not} include this section in the anonymized submission, only in the final paper. You can use the \texttt{ack} environment provided in the style file to automatically hide this section in the anonymized submission.

This research was supported by the Energy System Co-Design with Multiple Objectives and Power Electronics (E-COMP) Initiative, under the Laboratory  Directed Research and Development (LDRD) Program at Pacific Northwest National Laboratory (PNNL). %PNNL is a multi-program national laboratory operated for the U.S. Department of Energy (DOE) by Battelle Memorial Institute under Contract No. DE-AC05-76RL01830.
It was also supported by the Energy Earthshots Inititaive as part of the DOE Office of Biological and Environmental Research. PNNL is operated by DOE by the Battelle Memorial Institute under Contract DE-A06-76RLO 1830.

This work was also partially supported by the Ralph O’Connor Sustainable Energy Institute at Johns
Hopkins University, as well as NSF grants 2242931, 2232054, and 2143706.

\end{ack}

\newpage

\bibliographystyle{aaai2026}

\bibliography{aaai2026}


\newpage

\appendix

\section{Experimental Details: Particle Swarm Optimization}
\label{app:PSO_details}

Particle Swarm Optimization (PSO) is a metaheuristic optimization method which works by having a population (called a swarm) of candidate solutions (called particles). Particles update their solutions using simple update rules based on their own best-known position in search space, as well as the entire swarm's best-known position. PSO is commonly applied to optimization problems with complex objective functions and simple constraints. While simple bounds can be handled naturally in PSO, more complex constraints are often handled using penalty functions. A survey of constraint-handling techniques in PSO is found in \cite{innocente2021constraint}. 

For both control co-design experiments, we implement PSO baseline methods using the open-source library \textit{pyswarms} \cite{miranda2018pyswarms}.  We adopt the penalty-function approach to handle coupling constraints in our PSO baseline methods, and since the remaining upper-level constraints take the form of variable bounds, those are handled natively in the PSO algorithm of pyswarms. As is the case in the paper's main proposal, the lower-level problem must be feasible relative to the upper-level solutions found at each iteration of PSO. The constraints preventing this condition are treated as coupling constraints and enforced with penalty functions in the lower-level problem. For the control co-design problems, the lower-level problem implementations are identical to \eqref{appeq:twotank_UL} and  \eqref{appeq:building}. Overall, the PSO optimizes a relaxed upper-level objective function equal to 
\begin{equation}
    \label{eq:PSO_objective}
         \cL_{\bm{p}}(\bm{y}) + \kappa \nu(\bm{y}),  
\end{equation}
subject to $\bm{y} \in \cC_{\bm{p}}$, which are simple bounds on $\bm{y}$ in both of our experimental cases. Note that evaluation of $\nu(\bm{y})$ requires optimization at the lower level, at each step of PSO.

In each experiment, the pyswarms solver is given its default cognitive, social and inertia parameters $c_1 = 0.5, c_2 = 0.5, w = 0.9$, and run with $128$ particles for $200$ iterations.  The penalty coefficient $\kappa$ is chosen so that average coupling constraint violations over the test set are on the order of $1e-2$. This corresponds to $\kappa = 100.0$ in the two-tank experiment and $\kappa = 5.0$ in the HVAC experiment.

\section{Additional Results: Particle Swarm Optimization}
\label{app:PSO_results}

\begin{figure}[]
\centering
\includegraphics[width=0.66\columnwidth]{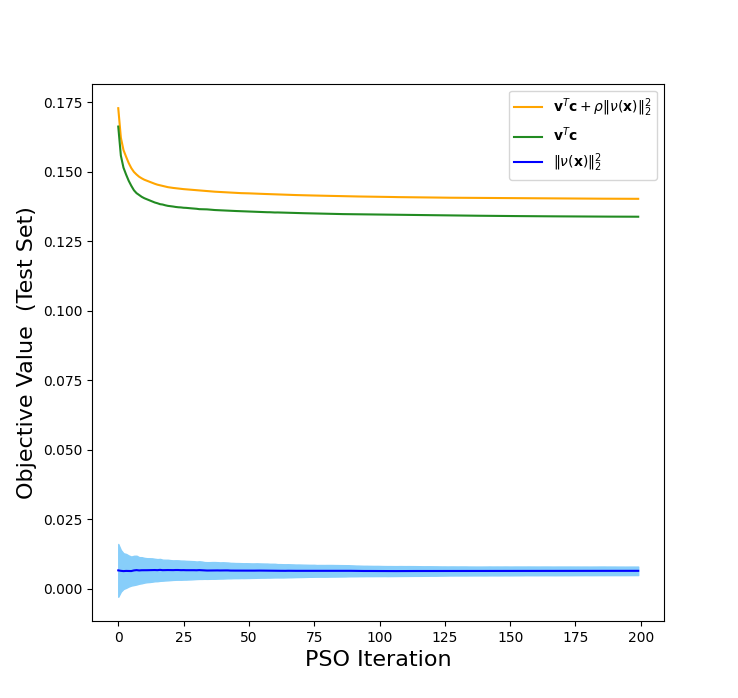}
\caption{Best Objective Value per PSO Iteration on Two-Tank System Co-Design }
\label{fig:tt_PSO_curve}
\end{figure}

\begin{figure}[]
\centering
\includegraphics[width=0.57\columnwidth]{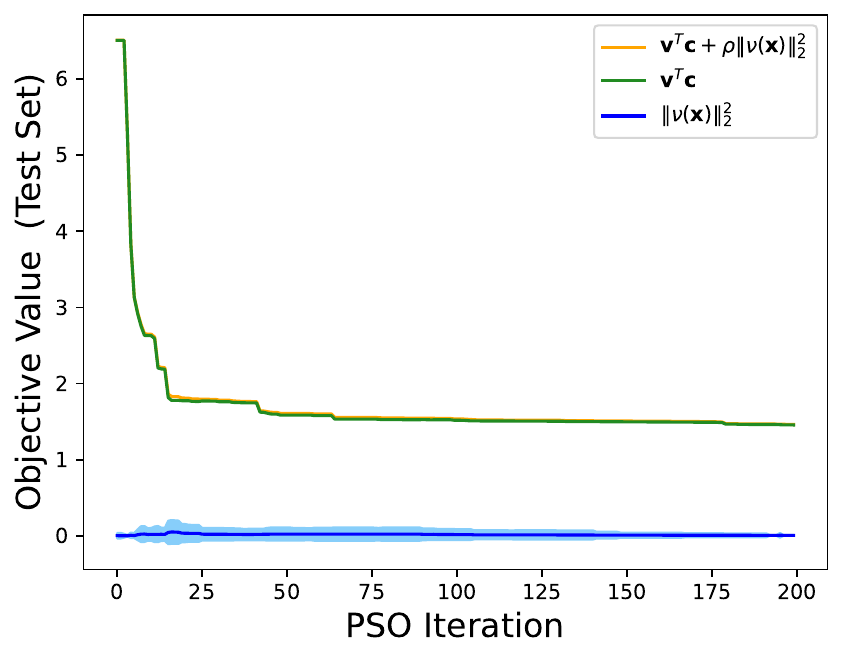}
\caption{Best Objective Value per PSO Iteration on HVAC System Co-Design }
\label{fig:PSO_curves_HVAC}
\end{figure}

We illustrate the evolution of the PSO objective throughout its solution of the test set instances. We plot the best objective values, among all particles, per iteration of PSO. The full PSO objective is illustrated, and also is shown in terms of upper-level design objective and coupling constraint violation penalty. Each metric is reported on average over the test set. Figure \ref{fig:tt_PSO_curve} corresponds to the two-tank problem, and \ref{fig:PSO_curves_HVAC} corresponds to the HVAC problem.

\section{Experimental Details: Learning Bilevel Quadratic Programming}
\label{app:BQP}

This section reports additional details on the experiments  presented in Section \ref{sec:bilevel_QP}.

\subsection{Hyperparameters and Training}
Results in Section \ref{sec:Experiments} are shown from the model which achieves the lowest loss among independent training runs using all combinations of the following hyperparameters: 

\begin{itemize}
    \item Learning rates from among $[ 10^{-1}, 10^{-2}, 10^{-3}, 10^{-4}, 10^{-5} ]$
    \item Correction stepsizes $\gamma$ from $[ 10^{-2}, 10^{-3}, 10^{-4}]$
    \item $\cL^{\texttt{SOFT}}$ penalty weights $\lambda$ from $[ 10^2, 10^3  ]$
\end{itemize}

The best values are $10^{-3}, 10^{-4}, 10^2$ respectively. All models are trained using the Adam optimizer \cite{zhang2018improved} in PyTorch. In each training run, $10$ correction steps are applied in training and $20$ are applied at test time.

\section{Experimental Details: Learning Control Co-design of a Two-Tank System}
\label{app:TT}

This section reports additional details on the nonlinear system design experiments  presented in Section \ref{sec:ControlCoDesign}.

\subsection{Problem Reformulation}

The two-tank system design and control problem (\ref{eq:twotank_UL}, \ref{eq:twotank_LL_discrete}) is bound by the coupling constraint $\bm{x}^{(N)} = \bm{p}$, which is redundantly placed at both levels to emphasize its coupling effect. This coupling constraint expresses that a valid system design must be controllable to the end-state $\bm{p}$. We recognize that this condition may not be satisfiable for any design variable $\bm{y}$; for instance, some $\bm{y}$ may not allow sufficient throughput to fill the tanks from $\bm{0}$ to $\bm{p}$ by time step $N$.

In practice, we therefore reformulate the problem as follows, as prefaced in Section \ref{sec:LowerLevelConstraints}:
\begin{minipage}{0.48\textwidth}
\begin{subequations}
\label{appeq:twotank_UL}
\begin{align}
    \mathcal{B}({\bm{p}}) = {\argmin}_{\bm{y}}\;\;\; \bm{v}^T \bm{y} \label{appeq:twotank_UL_obj} \\
     \textit{s.t.}  \;\;\; 
    \bm{x},\bm{u} = \mathcal{O}_{\bm{p}}(\bm{y})  \label{appeq:twotank_UL_opt} \\ 
    \;\;\;\;\;\;\;\;\;\;\bm{y}_{min} \leq \bm{y} \leq \bm{y}_{max} \label{appeq:twotank_bounds} \\
 \;\;\; \bm{x}^{(N)} = \bm{p}.  \label{appeq:twotank_UL_endpt}
\end{align}
\end{subequations}
\end{minipage}
\hfill
\begin{minipage}{0.48\textwidth}
\begin{subequations}
\label{appeq:twotank_LL_discrete}
\begin{align}
    \mathcal{O}_{\bm{p}}(\bm{y}) & \coloneqq \underset{0 \leq \bm{x}, \bm{u} \leq 1}{\argmin}   \;\;\;\sum_{k=1}^N \| \bm{u}^{(k)} \|_2^2 + \rho \|\bm{x}^{(N)} - \bm{p}\|^2 \\
     \textit{s.t.} \;\;\; & \bm{x}^{(i+1)} = \text{ODESolve}\big(f(\bm{x}^{(i)}, \bm{u}^{(i)}, \bm{y})\big) \label{appeq:twotank_LL_discrete_controls}
    %  & x_1^{(i+1)} = x_1^{(i)} + dt \cdot \left( y_1 (1-u_2^{(i)}) u_1^{(i)} - y_2 \sqrt{x_1^{(i)}} \right) \label{appeq:twotank_LL_discrete_states} \\
    % & x_2^{(i+1)} = x_2^{(i)} + dt \cdot \left(  y_1 u_2^{(i)} u_1^{(i)} + y_2 \left(\sqrt{x_1^{(i)}} -  \sqrt{x_2^{(i)}} \right) \right),  \label{appeq:twotank_LL_discrete_controls}
\end{align}
\end{subequations}
\end{minipage}
In this reformulation, the constraint $\bm{x}^{(N)} = \bm{p}$ remains at the upper level as a \emph{coupling} constraint, since it binds lower-level variables within the upper-level problem. Thus, it is treated by the coupling constraint correction in Algorithm \ref{alg:NCC}. The constraint is absent however, from the lower level in this formulation and instead replaced with a penalty on the lower-level objective. We take $\rho = 100$ as the penalty weight in all experimental runs. While equivalent to the original bilevel problem (\ref{eq:twotank_UL}, \ref{eq:twotank_LL_discrete}), the above formulation accomodates lower-level feasibility for any candidate design solution at the upper level.

\subsection{Hyperparameters and Training}

Results are shown from the model which achieves the lowest upper-level objective value on average, from among those whose average coupling constraint violation is less than or equal to that of the PSO baseline method after $10$ epochs. The results are chosen from among independent training runs using all combinations of the following hyperparameters: 

\begin{itemize}
    \item Learning rates from among $[ 10^{-1}, 10^{-2}, 10^{-3}, 10^{-4}, 10^{-5} ]$
    \item Correction stepsizes $\gamma$ from $[ 10^{-2}, 10^{-3}, 10^{-4}]$
    \item $\cL^{\texttt{SOFT}}$ penalty weights $\lambda$ from $[ 10, 10^2, 10^3  ]$
\end{itemize}
The best values are $10^{-3}, 10^{-2}, 10$ respectively. All models are trained using the Adam optimizer in PyTorch. In each training run, $5$ correction steps are applied in training and $10$ are applied at test time. The upper-level objective function has linear coefficients $\bm{v}$ consisting of all ones in each problem instance, meaning that the inlet and outlet valve coefficients should be minimized with equal priority.

\section{Experimental Details: Learning Control Co-design of a Building HVAC System}
\label{app:HVAC}

This section reports additional details on the HVAC system design experiments  presented in Section \ref{sec:ControlCoDesign}.

\subsection{Problem Reformulation}
The building HVAC design and control problem (\ref{eq:building_UL},\ref{eq:building}) is coupled by the thermal constraints
$\underline{\bm{p}}^{(k)} \leq \bm{w}^{(k)} \leq \overline{\bm{p}}^{(k)}$, which appear at both levels to emphasize their coupling role. Before Algorithm \ref{alg:NCC} can be applied, we recognize that those constraints may not be satisfiable when design variables $\bm{Y}$ prevent heat flows from converting to temperature changes rapidly enough to stay within those changing bounds. 

To arrive at an equivalent problem which ensures feasibility at the lower level for any $\bm{Y}$, we introduce slack variables $\underline{\bm{s}}^{(k)}$, $\overline{\bm{s}}^{(k)}$ to both sides of \eqref{eq:HVAC_LL_bds}  yielding

$$\underline{\bm{p}}^{(k)} - \underline{\bm{s}}^{(k)} \leq \bm{w}^{(k)} \leq \overline{\bm{p}}^{(k)} + \underline{\bm{s}}^{(k)},$$

along with a no-slack condition which maintains equivalence to the original problem:

$$ \underline{\bm{s}}^{(k)} = \overline{\bm{s}}^{(k)} = \bm{0} \;\;\;\;\; \forall k . $$

It is held at the upper level, and replaced in the lower level by a penalty term:

\begin{minipage}[t]{0.4\textwidth}
\centering
\begin{subequations}
    \label{appeq:building_UL}
    \begin{align}
         \cB(\bm{p}) = {\argmin} \;\; & \textit{Tr}(\bm{V}^T \bm{Y})   \\
         \textit{s.t.}\;\;\;\; & \bm{x},\bm{u},\bm{w} = \mathcal{O}_{\bm{p}}(\bm{Y}) \\
          & \bm{Y} \geq \bm{0} \\ 
          & \underline{\bm{s}}^{(k)} = \overline{\bm{s}}^{(k)} = \bm{0} \;\;\; \forall k \label{appeq:HVAC_UL_slacks}
    \end{align}
\end{subequations}
\label{appeq:UL}
\end{minipage}
\hfill
\begin{minipage}[t]{0.58\textwidth}
\centering
\begin{subequations}
    \label{appeq:building}
    \begin{align}
         \cO_{\bm{p}}(\bm{Y}) = \underset{\bm{x}, 0 \leq \bm{u} \leq 1, \bm{w}, \bm{s}}{\argmin} \;\; & \sum_{k \in \{1 \ldots N \}} \| \bm{u}^{(k)} \|_2^2 + \rho \left(  \sum_{k \in \{1 \ldots N \}} \| \underline{\bm{s}}^{(k)} \|_2^2 \;+  \sum_{k \in \{1 \ldots N \}} \| \overline{\bm{s}}^{(k)} \|_2^2  \right)  \\
         \textit{s.t.}\;\;\;\;  
          & \bm{w}^{(k)} = \bm{C}\bm{x}^{(k)} \\
          &  \underline{\bm{p}}^{(k)} - \underline{\bm{s}}^{(k)} \leq \bm{w}^{(k)} \leq \overline{\bm{p}}^{(k)} + \underline{\bm{s}}^{(k)}  \label{appeq:HVAC_LL_bds} \\
            &\bm{x}^{(k+1)} = \bm{A}\bm{x}^{(k)} + \bm{Y}\bm{u}^{(k)} + \bm{E}\bm{d}^{(k)} .
            \end{align}
\end{subequations}
\label{appeq:LL}
\vspace{-12pt}
\end{minipage}

In our implementation of Algorithm \ref{alg:NCC}, coupling constraint corrections are applied to \eqref{appeq:HVAC_UL_slacks}. It is the operative coupling constraint in this reformulation, binding the upper-level problem to lower-level variables $\underline{\bm{s}}^{(k)}, \overline{\bm{s}}^{(k)}$ .

\subsection{Hyperparameters and Training}
Results are shown from the model which achieves the lowest upper-level objective value on average, from among those whose average coupling constraint violation is less than or equal to that of the PSO baseline method after $25$ epochs. The results are chosen from among independent training runs using all combinations of the following hyperparameters: 

\begin{itemize}
    \item Learning rates from among $[ 10^{-1}, 10^{-2}, 10^{-3}, 10^{-4}, 10^{-5} ]$
    \item Correction stepsizes $\gamma$ from $[ 10^{-4}, 10^{-5}, 10^{-6}]$
    \item $\cL^{\texttt{SOFT}}$ penalty weights $\lambda$ from $[ 10^2, 10^3  ]$
\end{itemize}
The best values are $10^{-3}, 10^{-4}, 10^2$ respectively. All models are trained using the Adam optimizer in PyTorch. In each training run, $5$ correction steps are applied in training and $10$ are applied at test time. The upper-level objective function has linear coefficients $\bm{V}$ consisting of all ones in each problem instance, meaning that all elements of the actuator design variable have equal cost.

\end{document}